\documentclass[12pt]{amsart}
\usepackage{amssymb}
\usepackage{indentfirst,latexsym}
\usepackage{mathrsfs}
\usepackage{graphicx}
\usepackage{fancyhdr}
\usepackage{amsmath}

\date{}
\begin{document}


\title{The descent statistic over $123$-avoiding permutations}
\author{Marilena Barnabei, Flavio Bonetti, and Matteo Silimbani}
\address{
 Department of Mathematics, University of Bologna\\ P.zza di Porta San Donato 5, 40126 Bologna}
 \email{barnabei@dm.unibo.it\\bonetti@dm.unibo.it\\silimban@dm.unibo.it}\maketitle

 \noindent {\bf Abstract} We exploit Krattenthaler's bijection between $123$-avoiding permutations and Dyck paths
 to determine the Eulerian distribution over the set $S_n(123)$ of $123$-avoiding permutations in $S_n$. In particular, we show that the descents of a permutation
 correspond to valleys and triple falls of the associated Dyck path. We get the Eulerian numbers of $S_n(123)$ by studying the joint
 distribution of these two statistics on Dyck paths.
\newline

\newline


\section{Introduction}

\noindent A permutation $\sigma\in S_n$ \emph{avoids a pattern}
$\tau\in S_k$ if $\sigma$ does not contain a subsequence that is
order-isomorphic to $\tau$. The subset of $S_n$ of all
permutations avoiding a pattern $\tau$ is denoted by $S_n(\tau)$.
Pattern avoiding permutations have been intensively studied in
recent years from many points of view (see e.g. \cite{kratt},
\cite{simisch} and references therein).

\noindent In the case $\tau\in S_3$, it has been shown that the
cardinality of $S_n(\tau)$ equals the $n$-th Catalan number, for
every pattern $\tau$, and hence the set $S_n(\tau)$ is in
bijection with the set of Dyck paths of semilength $n$. Indeed,
the six patterns in $S_3$ are related as follows:
\begin{itemize}
\item $321=123^{rev}$,
\item $231=132^{rev}$,
\item $213=132^c$,
\item $312=(132^{c})^{rev}$,
\end{itemize}
where $rev$ and $c$ denote the usual reverse and complement
operations. Hence, in order to determine the distribution of the
descent statistic over $S_n(\tau)$, for every $\tau\in S_3$, it is
sufficient to examine the distribution of descents over two sets
$S_n(132)$ and $S_n(123)$.

\noindent In both cases, the two bijections due to Krattenthaler
\cite{kratt} allow to translate the descent statistic into some
appropriate statistics on Dyck paths.

\noindent In the case $\tau=132$, the descents of a permutation
are in one-to-one correspondence with the valleys of the
associated Dyck path (see \cite{oeis}).\\

\noindent In this paper we investigate the case $\tau=123$. In
particular, we exploit a variation of Krattenthaler's map to
translate the descents of a permutation $\sigma\in S_n(123)$ into
peculiar subconfigurations of the associated Dyck path, namely,
valleys and triple falls.

\noindent For that reason, we study the joint distribution of
valleys and triple falls over the set $\mathcal{P}_n$ of Dyck
paths of semilength $n$, and we give an explicit expression for
its trivariate generating function
$$A(x,y,z)=\sum_{n\geq 0}\sum_{\mathscr{D}\in \mathcal{P}_n}x^ny^{v(\mathscr{D})}z^{tf(\mathscr{D})}=\sum_{n,p,q\geq 0} a_{n,p,q}x^ny^pz^q,$$
where $v(\mathscr{D})$ denotes the number of valleys in
$\mathscr{D}$ and $tf(\mathscr{D})$ denotes the number of triple
falls in $\mathscr{D}$. This series specializes into some well
known generating functions, such as the generating function of
Catalan numbers, Motzkin numbers, Narayana numbers, and seq.
A$092107$ in \cite{oeis} (see also \cite{sapou}).

\section{Dyck paths}

\noindent A \emph{Dyck path} is a lattice path in the integer
lattice $\mathbb{N}\times\mathbb{N}$ starting from the origin,
consisting of up-steps $U=(1,1)$ and down steps $D=(1,-1)$, never
passing below the x-axis, and ending at ground level.

\noindent We recall that a \emph{return} of a Dyck path is a down
step ending on the $x$-axis. An \emph{irreducible}  Dyck path is a
Dyck path
 with exactly one return.\\

\noindent We observe that a Dyck path $\mathscr{D}$ can be
decomposed according to its last return (\emph{last return
decomposition}) into the juxtaposition of a (possibly empty) Dyck
path $\mathscr{D}'$ of shorter length and an irreducible Dyck path $\mathscr{D}''$.\\

\begin{figure}[ht]
\begin{center}
\includegraphics[bb=58 569 471 716,width=.8\textwidth]{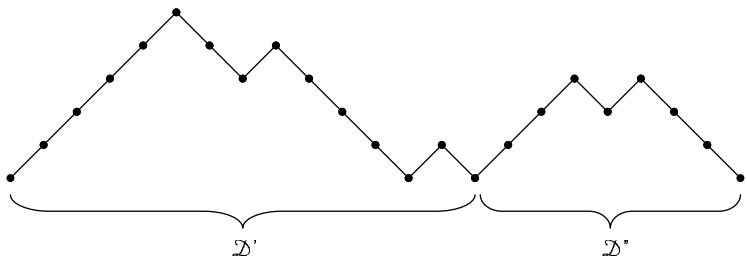} \caption{The
last return decomposition of the Dyck path
$\mathscr{D}=U^5D^2UD^4UDU^3DUD^3$.}\label{fighiur}
\end{center}
\end{figure}

\noindent For example, the Dyck path
$\mathscr{D}=U^5D^2UD^4UDU^3DUD^3$ decomposes into
$\mathscr{D}'\bigoplus\mathscr{D}''$, where
$\mathscr{D}'=U^5D^2UD^4UD$ and $\mathscr{D}''=U^3DUD^3$, as shown
in Figure \ref{fighiur}.

\section{Krattenthaler's bijection}

\noindent In \cite{kratt}, Krattenthaler describes a bijection
between the set $S_n(123)$ and the set $\mathcal{P}_n$ of Dyck
paths of semilength $n$. We present a slightly modified version of
this bijection.

\noindent Let $\sigma=\sigma(1)\ldots\sigma(n)$ be a
$123$-avoiding permutation. Recall that a \emph{left-to-right
minimum} of $\sigma$ is an element $\sigma(i)$ which is smaller
than $\sigma(j)$, with $j<i$ (note that the first entry
$\sigma(1)$ is a left-to-right minimum). Let $x_1,\ldots, x_s$ be
the left-to-right minima in $\sigma$. Then, we can write
\begin{equation}\sigma=x_1\,w_1\,\ldots\,x_s\,w_s,\label{comek}\end{equation}
where $w_i$ are (possibly empty) words. Moreover, since $\sigma$
 avoids $123$, the word $w_1\,w_2\,\ldots\,w_s$ must be
 decreasing.\\

 \noindent In order to construct the Dyck path $\kappa(\sigma)$ corresponding to $\sigma$, read the decomposition (\ref{comek}) from left to
 right. Any left-to-right minimum $x_i$ is translated into
 $x_{i-1}-x_i$ up steps (with the convention $x_0=n+1$) and any
 subword $w_i$ is translated into $l_i+1$ down steps, where $l_i$
 denotes the number of elements in $w_i$.

\begin{figure}[ht]
\begin{center}
\includegraphics[bb=61 602 324 695,width=.7\textwidth]{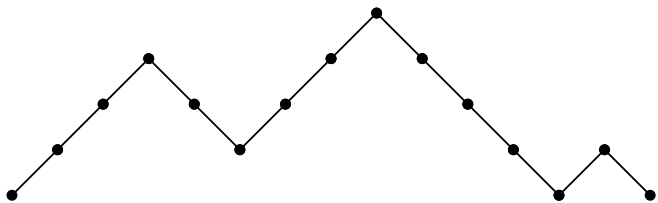} \caption{The
Dyck path $\kappa(\sigma)$, with
$\sigma=5\,7\,2\,6\,4\,3\,1$.}\label{ballenji}
\end{center}
\end{figure}

 \noindent For example, the permutation $\sigma=5\,7\,2\,6\,4\,3\,1$ in $S_7(123)$ corresponds to
 the path in Figure \ref{ballenji}.

\section{The descent statistic}


\noindent We say that a permutation $\sigma$ has a \emph{descent}
at position $i$ if $\sigma(i)>\sigma(i+1)$. We denote by
des$(\sigma)$ the number of descents of the permutation
$\sigma$.\\

\noindent In this section we determine the generating function
$$E(x,y)=\sum_{n\geq 0}\sum_{\sigma\in S_n(123)}x^ny^{\textrm{des}(\sigma)}=\sum_{n\geq 0}\sum_{k\geq 0}e_{n,k}x^ny^k,$$
where $e_{n,k}$ denotes the number of permutations in $S_n(123)$
with $k$ descents.

 \newtheorem{inizio}{Proposition}
\begin{inizio}\label{sfrutt}
Let $\sigma$ be a permutation in $S_n(123)$, and
$\mathscr{D}=\kappa(\sigma)$. The number of descents of $\sigma$
is
$$\textrm{des}(\sigma)=v(\mathscr{D})+tf(\mathscr{D}),$$
where $v(\mathscr{D})$ is the number of valleys in $\mathscr{D}$
and $tf(\mathscr{D})$ is the number of \emph{triple falls} in
$\mathscr{D}$, namely, the number of occurrences of $DDD$ in
$\mathscr{D}$.\end{inizio}

\noindent \emph{Proof} Let $\sigma=x_1\,w_1\,\ldots\,x_s\, w_s$ be
a $123$-avoiding permutation.  The descents of $\sigma$ occur
precisely in the following positions:
\begin{itemize}
\item[1.] between two consecutive symbols in the same word
$w_i$ (we have $l_i-1$ of such descents),
\item[2.] before every left-to-right
minimum $x_i$, except for the first one.
\end{itemize}
The proof is completed as soon as we remark that:
\begin{itemize}
\item[1.] every word $w_i$ is mapped into a descending run of
$\kappa(\sigma)$ of length $l_i+1$. Such descending run contains
$l_i-1$ triple falls, that are therefore in bijection with the
descents contained in $w_i$,
\item[2.] every left-to-right
minimum $x_i$ with $i\geq 2$ corresponds to a valley in
$\kappa(\sigma)$.
\end{itemize}
\begin{flushright}
$\diamond$
\end{flushright}

\noindent The preceding result implies that we can switch our
attention from permutations in $S_n(123)$ with $k$ descents to
Dyck paths of semilength $n$ with $k$ among valleys and triple
falls. Hence, we study the joint distribution of valleys and
triple falls over $\mathcal{P}_n$, namely, we analyze the
generating function
$$A(x,y,z)=\sum_{n\geq 0}\sum_{\mathscr{D}\in \mathcal{P}_n}x^ny^{v(\mathscr{D})}z^{tf(\mathscr{D})}=\sum_{n,p,q\geq 0} a_{n,p,q}x^ny^pz^q.$$
We determine the relation between the function $A(x,y,z)$ and the
generating function
$$B(x,y,z)=\sum_{n\geq 0}\sum_{\mathscr{D}\in \mathcal{IP}_n}x^ny^{v(\mathscr{D})}z^{tf(\mathscr{D})}=\sum_{n,p,q\geq 0} b_{n,p,q}x^ny^pz^q$$
of the same joint distribution over the set $\mathcal{IP}_n$ of
irreducible Dyck paths in $\mathcal{P}_n$.

\newtheorem{manfri}[inizio]{Proposition}
\begin{manfri}\label{struff}
For every $n>2$, we have:
\begin{equation}b_{n,p,q}=a_{n-1,p,q-1}-a_{n-2,p-1,q-1}+a_{n-2,p-1,q}.\label{bagno}\end{equation}
\end{manfri}

\noindent \emph{Proof} An irreducible Dyck path of semilength $n$
with $p$ valleys and $q$ triple falls can be obtained by
prepending $U$ and appending $D$ to a Dyck path of semilength
$n-1$ of one of the two following types:
\begin{itemize}
\item[1.] a Dyck path with $p$ valleys and $q$ triple falls, ending
with $UD$,
\item[2.] a Dyck path with $p$ valleys and $q-1$ triple falls, not
ending with $UD$.
\end{itemize}
We remark that:
\begin{itemize}
\item[1.] the paths of the first kind are in bijection with Dyck
paths of semilength $n-2$ with $p-1$ valleys and $q$ triple falls,
enumerated by $a_{n-2,p-1,q}$.
\begin{figure}[h]
\begin{center}
\includegraphics[bb=61 623 379 739,width=.8\textwidth]{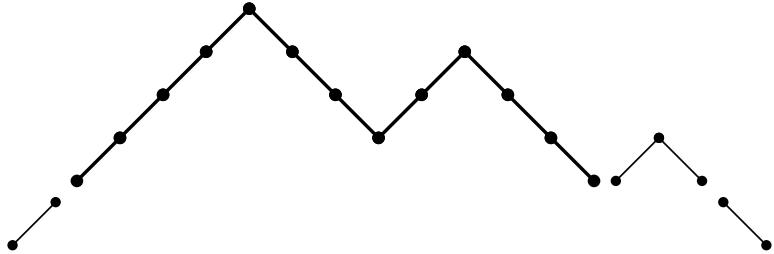} \caption{The Dyck path $U^5D^3U^2D^3UD^2$ with $2$ valleys and $2$ triple falls is obtained by
appending $UD$ to the path $U^4D^3U^2D^3$ with $1$ valley and $2$
triple falls, and then elevating.}
\end{center}
\end{figure}
\item[2.] in order to enumerate the paths of the second kind we
have to subtract from the integer $a_{n-1,p,q-1}$ the number of
Dyck paths of semilength $n-1$ with $p$ valleys and $q-1$ triple
falls, ending with $UD$. Dyck paths of this kind are in bijection
with Dyck paths of semilength $n-2$ with $p-1$ valleys and $q-1$
triple falls, enumerated by $a_{n-2,p-1,q-1}$.
\begin{figure}[h]
\begin{center}
\includegraphics[bb=61 623 379 739,width=.8\textwidth]{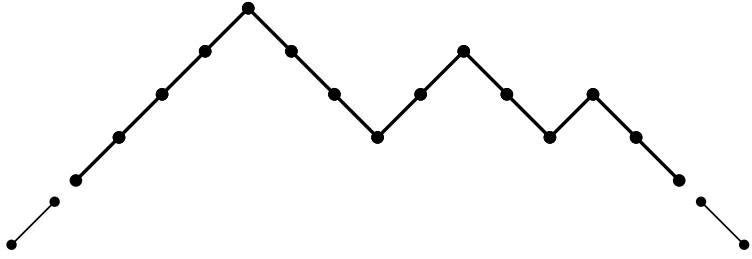} \caption{The Dyck path $U^5D^3U^2D^2UD^3$ with $2$ valleys and $2$ triple falls
is obtained by elevating the path $U^4D^3U^2D^2UD^2$ with $2$
valleys and $1$ triple fall.}
\end{center}
\end{figure}
\end{itemize}
\begin{flushright}
$\diamond$
\end{flushright}

\newtheorem{rusco}[inizio]{Proposition}
\begin{rusco}\label{ripe}
For every $n>0$, we have:
\begin{equation}a_{n,p,q}=b_{n,p,q}+\sum_{i=1}^{n-1}\sum_{j,s>0}b_{i,j,s}a_{n-i,p-j-1,q-s}.\label{vabeh}\end{equation}
\end{rusco}

\noindent \emph{Proof} Let $\mathscr{D}$ be a Dyck path of
semilength $n$ and consider its last return decomposition
$\mathscr{D}=\mathscr{D}'\bigoplus\mathscr{D}''$. If
$\mathscr{D}'$ is empty, then $\mathscr{D}$ is irreducible.
Otherwise:
\begin{itemize}
\item $v(\mathscr{D})=v(\mathscr{D}')+v(\mathscr{D}'')+1$,
\item $tf(\mathscr{D})=tf(\mathscr{D}')+tf(\mathscr{D}'')$.
\end{itemize}
\begin{flushright}
$\diamond$
\end{flushright}

\noindent Identities (\ref{bagno}) and (\ref{vabeh}) yield the
following relations between the two generating functions
$A(x,y,z)$ and $B(x,y,z)$:

\newtheorem{conc}[inizio]{Proposition}
\begin{conc}\label{arch}
We have:
\begin{equation}
B(x,y,z)=(A(x,y,z)-1)(xz+x^2y-x^2yz)+1+x+x^2-x^2z,\label{uno}\end{equation}
\begin{equation}A(x,y,z)=B(x,y,z)+y(B(x,y,z)-1)(A(x,y,z)-1)\label{due}\end{equation}
\end{conc}

\noindent \emph{Proof} Observe that recurrence (\ref{bagno}) holds
for $n>2$. This fact gives rise to the correction terms of degree
less than $3$ in Formula (\ref{uno}).
\begin{flushright}
$\diamond$
\end{flushright}

\noindent Combining Formul\ae (\ref{uno}) and (\ref{due}) we
obtain the following:

\newtheorem{princ}[inizio]{Theorem}
\begin{princ}\label{arch}
We have: {\setlength\arraycolsep{2pt}
\begin{eqnarray}
A(x,y,z) & = & \frac{1}{2xy(xyz-z-xy)}\left(-1+xy+2x^2y\right.
\nonumber\\
& &-2x^2y^2+xz-2xyz-2x^2yz+2x^2y^2z\\
& &
\left.+\sqrt{1-2xy-4x^2y+x^2y^2-2xz+2x^2yz+x^2z^2}\right)\nonumber
\end{eqnarray}}
\end{princ}
\begin{flushright}
$\diamond$
\end{flushright}

\noindent This last result allows us to determine the generating
function $E(x,y)$ of the Eulerian distribution over $S_n(123)$. In
fact, previous arguments show that $$E(x,y)=A(x,y,y)$$ and hence:

\newtheorem{scopo}[inizio]{Theorem}
\begin{scopo}\label{arch}
We have:
$$E(x,y)=\frac{-1+2xy+2x^2y-2xy^2-4x^2y^2+2x^2y^3+\sqrt{1-4xy-4x^2y+4x^2y^2}}{2xy^2(xy-1-x)}.$$
\end{scopo}
\begin{flushright}
$\diamond$
\end{flushright}

\noindent The first values of the sequence $e_{n,d}$ are shown in
the following table:

$$\begin{array}{l|lllllll}
n/d&0&1&2&3&4&5&6\\\hline
0&1&&&&&&\\
1&1&&&&&&\\
2&1&1&&&&&\\
3&0&4&1&&&&\\
4&0&2&11&1&&&\\
5&0&0&15&26&1&&\\
6&0&0&5&69&57&1&\\
7&0&0&0&56&252&120&1\\\vspace{1cm}
\end{array}$$

\noindent Needless to say, the series $A(x,y,z)$ specializes into
some well known generating functions. In particular, $A(x,1,1)$ is
the generating function of Catalan numbers, $A(x,1,0)$ the
generating function of Motzkin numbers, $yA(x,y,1)$ the generating
function of Narayana numbers, and $A(x,1,z)$ the generating
function of seq. A$092107$ in \cite{oeis}.

\end{document}